\def\IR{{\Bbb R}}

\def\IC{\Bbb C} 
\def\ID{{\Bbb D}}

\def\zbar{{\overline{z}}} 
\def\zetabar{{\overline{\zeta}}} 
\def\wbar{{\overline{w}}}

\documentclass[12pt]{article} 

\usepackage[psamsfonts]{amssymb} 
\usepackage{amsfonts,amsmath}

\usepackage{epsfig,multicol}

\newtheorem{theorem}{Theorem}
\newtheorem{corollary}{Corollary}

\title{Super regularity for Beltrami systems.}
\author{ Gaven J. Martin  } 
\date{} 
\begin{document}

\maketitle

\begin{abstract}  
We prove a surprising higher regularity for solutions to the nonlinear elliptic autonomous Beltrami equation in a planar domain $\Omega$,
\[ f_\zbar = {\cal A}(f_z) \hskip15pt a.e.\;\; z\in \Omega, \]
when ${\cal A}$ is linear at $\infty$.  Namely $W^{1,1}_{loc}(\Omega)$ solutions are $W^{2,2+\epsilon}_{loc}(\Omega)$.  Here $\epsilon>0$ depends explicitly on the ellipticity bounds of ${\cal A}$.  The condition ``is linear at $\infty$'' is necessary - the result is false for the equation $f_\zbar = k|f_z|$,  for any $0<k<1$, ($k=0$ is Weyl's lemma).

We discuss the subsequent higher regularity implications for fully non-linear Beltrami systems.
\end{abstract}

\section{Introduction}     

The governing equations of planar geometric analysis and the theory of quasiconformal mappings,  Teichm\"uller spaces and so forth are the Beltrami equations and their nonlinear counterparts,  see for instance  \cite{AIM,IM,G1,G2,H,HM} and elsewhere.  Beltrami equations come in several different flavours.  As examples,  let $\Omega\subset \IC$ be a domain and let $f:\Omega\to\IC$ be a mapping of Sobolev class $W^{1,1}_{loc}(\Omega)$ consisting of functions whose first derivatives are locally integrable,  then we have
\begin{itemize}
\item $\IC$-linear:  $f_\zbar=\mu(z) f_z$,  with ellipticity estimate $\|\mu\|_{L^\infty(\Omega)}<1$;
\item $\IR$-linear:  $f_\zbar=\mu(z) f_z+\nu(z) \overline{f_z}$,  with ellipticity estimate \[ \|\,|\mu|+|\nu|\,\|_{L^\infty(\Omega)}<1;\]
\item Autonomous: $f_\zbar={\cal A}(f_z)$,  with ellipticity estimate: there is $k<1$ so that for all $\zeta,\eta\in\IC$ \[|{\cal A}(\zeta)-{\cal A}(\eta)|\leq k |\zeta-\eta|\]
\item Fully nonlinear: $f_\zbar={\cal H}(z,f,f_z)$,  with ellipticity estimate : there is $k<1$ so that for all $z\in \Omega$,  all $\zeta,\eta\in\IC$ \[|{\cal H}(z,w,\zeta)-{\cal H}(z,w,\eta)|\leq k |\zeta-\eta|\]
with additional conditions on ${\cal H}$,  see \cite[Chapters 7 \& 8]{AIM}.
\end{itemize}
A homeomorphic $W^{1,2}_{loc}(\Omega)$ solution to any such an equation is a {\em quasiconformal mapping}.  The theory of quasiconformal mappings is comprehensively treated in \cite{AIM}. An orientation-preserving  homeomorphism   $f:\Omega \to
\Omega^\prime$ is  $K$-{\em quasiconformal},  $1  \leqslant K<\infty$,  if  $f \in W^{1,2}_{loc}(\Omega)$
 and 
\begin{equation} 
\max_{\alpha} |\partial_\alpha f(z)|   \leqslant K \min_{\alpha} |\partial_\alpha f(z)|
\end{equation}
for almost every $z\in \Omega$.  

\medskip
 
 More generally,  any $W^{1,2}_{loc}(\Omega)$ solution $f$ to any of the above equations is {\em quasiregular} and factors as $f=\Phi\circ g$ with $g$ quasiconformal and $\Phi:g(\Omega)\to \IC$ conformal - the Sto\"{i}low factorisation theorem, \cite{St}.
 
Each of these equations has a seminal application and they are all inter-related.  The apriori assumption that $f\in W^{1,1}_{loc}(\Omega)$ is so that we can even speak of $f$ as a ``solution''.  Without stronger assumptions on $\mu$ or ${\cal H}$ not much can be said,  but note for instance that $\mu=0$ on an open set implies $f$ is holomorphic on that set - Weyl's Lemma.  The higher regularity theory of these equations typically assumes more on $f$,  for instance $f\in W^{1,q}_{loc}(\Omega)$ for some $1<q \leq 2 $ usually depending on the ellipticity constant $k$,  and in return delivers a far nicer outcome,  $f\in W^{1,p}_{loc}(\Omega)$ for some $p > 2$,  again depending on $k$.  Astala's theorem \cite{Ast} gives the optimal result in the $\IC$-linear case and can be used to analyse other cases.  Questions of existence and uniqueness are fairly well understood through the topological properties of these mappings and Sto\"{i}low factorisation,  see also \cite[\S 5.5 \& \S 6.1]{AIM}.  However there are intriguing subtleties in the nonlinear case,  see \cite{AC1,AC2}.

\medskip
Here we present a rather surprising higher regularity theorem for autonomous systems.  We say that ${\cal A}:\IC\to \IC$ is linear at infinity if there are constants $a,b\in \IC$,  $|a|+|b|<1$,  and $\alpha$ with $0\leq \alpha<1$,  such that
\begin{equation}
{\cal A}(\zeta) = a \zeta + b \bar \zeta + O(|\zeta|^\alpha).
\end{equation}
The next theorem makes no ellipticity assumptions on ${\cal A}$.
\begin{theorem}[Super-regularity for Autonomous systems]\label{theorem1}
Let $\Omega\subset \IC$ be a planar domain and $f:\Omega\to \IC$ be a $W^{1,1}_{loc}(\Omega)$ solution to the autonomous Beltrami system
\begin{equation}
f_\zbar = {\cal A}(f_z),  \hskip20pt a.e. \; z\in \Omega \label{eqn2}
\end{equation}
where ${\cal A}$ is linear at $\infty$.  Then $f\in W^{1,p}_{loc}(\Omega)$ for all $p<\infty$.  

\medskip

If in addition $|{\cal A}(\zeta)|\leq k|\zeta|$ for some $k<1$ and all $\zeta\in\IC$,  then $f$ is also quasiregular.
\end{theorem}
In the next corollary we should realise that there is no assumed connection between $|a|+|b|$ and $k$.  However in applications it be obvious  that $|a|+|b|\leq k$.
\begin{corollary} \label{cor1} Suppose that $f\in W^{1,1}_{loc}(\Omega)$ is a solution to (\ref{eqn2}) for ${\cal A}$ linear at $\infty$,  and $f$ satisfies a Lipschitz bound of the form  
\begin{equation}\label{eqn3}
|f_\zbar(z+t\zeta)-f_\zbar(z) | \leq k \,  |f_z(z+t\zeta)-f_z(z) | , \hskip15pt |\zeta|=1,
\end{equation}
for some $k<1$ and for all $0< t< a(z)$ for some continuous function $a:\Omega\to \IR_+$,   $a(z)\leq {\rm dist}(z,\partial \Omega)$.
Then 
\begin{enumerate}
\item $f\in W^{2,q}_{loc}(\Omega)$ for all $q<1+1/k$,  
\item Each member of the $\IR$-linear family
\[ \{ a f_x(z)+b f_y(z) : a,b\in \IR \} \]
is either $\frac{1+k}{1-k}$-quasiregular mapping or a constant.
\item There are measurable $\mu,\nu:\Omega\to \IC$ with $|\mu|+|\nu|\leq k$ so that both directional derivatives $f_x$ and $f_y$ satisfy the $\IR$-linear Beltrami equation,
\[ h_\zbar = \mu(z) h_z +\nu(z) \overline{h_z}, \hskip20pt h\in \{f_x,f_y\}  \]
\item The complex gradient $f_z$ is itself quasiregular and satisifies the $\IR$ linear equation
\[ h_\zbar = \frac{\mu(z)}{1-|\nu(z)|^2}  h_z + \frac{\overline{\mu(z)}\nu(z)}{1-|\nu(z)|^2}  \overline{h_z}, \hskip20pt h = f_z  \]
and thus $f_z \in W^{1,q}_{loc}(\Omega)$ for all 
\[ q < s = 1+1/k',  \hskip15pt k'= \left\| \frac{|\mu|}{1-|\nu|} \right\|_\infty \leq k \]
\end{enumerate}
 \end{corollary}
 
 Of course the Lipschitz bound at (\ref{eqn3}) is achieved if for instance ${\cal A}$ satisfies the usual ellipticity assumption: ${\cal A}$ is $k$-Lipschitz,
 \begin{equation}
| {\cal A}(\zeta)-{\cal A}(\eta)|\leq k |\zeta-\eta|, \hskip20pt \zeta,\eta\in \IC.
 \end{equation}
 In this case one can make further claims about existence, uniqueness and higher regularity of $W^{1,2}_{loc}(\Omega)$ solutions \cite{AIM, AC1,AC2,JJ}.  An example of these which follows from the Schauder theory would be the following.
 
 \begin{corollary} \label{cor2} Let ${\cal A}:\IC\to\IC$ be smooth and linear at $\infty$ with 
 \begin{equation}
 |A_\zeta(\eta)|+|A_{\bar\zeta}(\eta)|<1, \hskip15pt \eta\in \IC.
 \end{equation}
 Then every $W^{1,1}_{loc}(\Omega)$ solution $f$ to the equation (\ref{eqn2}) is smooth and quasiregular.  After appropriate normalisation a homeomorphic $f$ will be unique.
 \end{corollary}
 
 In contrast to this corollary however is another very interesting result of Astala et al \cite{Astalaetal} which shows the condition of being linear at $\infty$ is essential.
 
 \begin{theorem} For each $k<1$ is an $f\in W^{1,1}_{loc}(\IC)$ solving the equation
\begin{equation} f_\zbar = k|f_z| \label{eqn6} \end{equation}
 for which $f\not\in W^{1,1+k}_{loc}(\IC)$.  Any solution to (\ref{eqn6}) which lies in $W^{1,q}_{loc}(\IC)$ for some $q>1+k$ is smooth.
 \end{theorem}

\section{Proofs.}

The main result,  Theorem \ref{theorem1},  and its corollaries are a consequence of an induction based on the following result.

\begin{theorem}\label{theorem3}
Let ${\cal A}:\IC\to \IC$ be linear at $\infty$ and $h\in L^{q}_{loc}(\Omega)$.  Then every $W^{1,1}_{loc}(\Omega)$ solution to the equation
\begin{equation}\label{eqn8}
f_\zbar = {\cal A}(f_z)+ h
\end{equation}
lies in $W^{1,q}_{loc}(\Omega)$.
\end{theorem}
 \noindent{\bf Proof.}  Let $\eta\in C_{0}^{\infty}(\Omega)$.  We write (\ref{eqn8})  as 
\[  f_\zbar = a f_z + b \overline{f_z} + O(|f_z|^\alpha) + h \]
 for some $\alpha<1$.  We multiply this equation by $\eta$ and use the fact that
 \[ (\eta f)_\zbar = \eta f_\zbar+ \eta_\zbar f, \hskip10pt   (\eta f)_z = \eta f_z+ \eta_z f \]
 and rearrange terms to achieve the following equation for $\tilde{f}=\eta f$.
 \begin{equation}\label{eqn9}
\tilde{f}_\zbar = a \tilde{f}_z + b \overline{\tilde{f}_z} + u,
 \end{equation}
 where 
 \[  u = \eta h + \eta O(|f_z|^\alpha) + \eta_\zbar f + a \eta_z f + b \overline{\eta_z f}  \]
 The Sobolev embedding gives $\eta f\in L^r(\IC)$ for every $r<2$ and so we can assume $u \in L^{1/\alpha}(\IC) \cap L^q(\IC) \cap L^r(\IC)$.  We therefore have to show that any compactly supported $W^{1,1}(\IC)$ solution to the inhomogeneous {\em constant coefficient} equation (\ref{eqn9}) has improved regularity.  If ${\cal S}$ denotes the Beurling transform, a singular integral operator of Calderon-Zygmund type,  since ${\cal S}\circ \frac{\partial}{\partial \zbar} = \frac{\partial}{\partial z}$ we really want to establish the invertibility of the constant coefficient Beltrami operator $I-a{\cal S} - b\overline{\cal S}$ since (\ref{eqn9}) reads as 
 \[ I-a {\cal S}(f_\zbar)- b \overline{ {\cal S}(f_\zbar)} = u \]
 More on Beltrami operators can be found in \cite[Chapter 4]{AIM}.  
 
 \medskip
 
 Following \cite[\S 15.2]{AIM} we address (\ref{eqn9})   by a linear change of variables reducing it to  the inhomogeneous Cauchy-Riemann equation.  
Namely, for any given constants $\mu, \nu \in \ID$ we  may use the transformation 
\[
\tilde{f}(\zeta) =  \frac{g(z)-\nu \overline{g(z)}}{1-|\nu|^2}, \quad \quad \zeta = z + \mu \zbar, 
\]
which rearranges to
 \begin{equation}
\label{afterahl22}
g(z)=\tilde{f}(\zeta)+\nu\overline{\tilde{f}(\zeta)},   \quad \quad z=\frac{\zeta-\mu\zetabar}{1-|\mu|^2}
\end{equation}
and we set $v(z)=u(\zeta) =u(z+\mu\zbar)$.   If we make the following choices for $\mu$ and $\nu$, 
 \begin{eqnarray*}
\mu & = & \frac{-2a}{1+|a|^2-|b|^2+\sqrt{(1+|a|-|b|^2)^2-4|a|^2}}\\ 
\nu & = & \frac{-2b}{1+|b|^2-|a|^2+\sqrt{(1+|b|-|a|^2)^2-4|a|^2}},
 \end{eqnarray*}
then since $|a|+|b|<1$ some computation reveals that indeed
\[ |\mu|  \leqslant \frac{|a|}{1-|b|^2} < 1,  \quad  \quad |\nu|  \leqslant \frac{|b|}{1-|a|^2}  <1. \]
Thus we can make these transformations above and  a few further elementary computations shows that  (\ref{eqn9}) now reads as 
\begin{equation} \frac{\partial g}{\partial \zbar}(z)=v(z)+a b\; \overline{v(z)} \end{equation} 
Hence $g_\zbar\in L^{1/\alpha}(\IC) \cap L^q(\IC) \cap L^r(\IC)$ and since ${\cal S}:L^q(\IC)\to L^q(\IC)$ is bounded,  we have $g_z={\cal S}(g_\zbar)\in L^{1/\alpha}(\IC) \cap L^q(\IC) \cap L^r(\IC)$.  Thus $g$ (and hence $\tilde{f}$ and consequently $f$) lies in the space $W^{1,s}_{loc}(\Omega)$,  where $s=\min\{1/\alpha, q, r\}$.  Iterating this construction we achieve $s=q$ as $\alpha\to \alpha^2$ and $r$ increases to the new Sobolev embedding exponent. \hfill $\Box$

\medskip

With this in hand Theorem \ref{theorem1} is immediate.  Corollary \ref{cor1} also follows as Theorem \ref{theorem1} promotes $f$ to an element of $W^{1,2}_{loc}(\Omega)$ (actually better) and so we may appeal to the main result of  \cite{HM}.  As for Corollary \ref{cor2},  the only remark that needs to be made is that if ${\cal A}$ is smooth with the proposed bounds, and linear at $\infty$,  then ${\cal A}$ must satisfy some Lipschitz bound with constant less than $1$.  Corollary \ref{cor2} then put the solution in $W^{2,2}_{loc}(\Omega)$,  the complex gradient is quasiregular,  and so H\"older continuous.  Then, as noted,  the Schauder bounds \cite{AIM} give the desired conclusion. 

\medskip

\section{A remark on the Hodographic transform.}  It is interesting to note that under the Hodographic transformation \cite[\S 16.3]{AIM} (basically a change of variables for homeomorphic solutions) and under a few regularity assumptions we see a relationship between conditions on  $f$ and $h=f^{-1}$.   The condition ``linear at $\infty$'' is simply that there is $q>1$ such that $|f_\zbar -a f_z-b \overline{f_z}|\in L^{q}_{loc}(\Omega)$.  If $h$ is sufficiently regular so as to change variables we find
\[ |f_\zbar(h) -a f_z(h)-b \overline{f_z(h)}| J(w,h) = |h_\wbar -b h_w- a \overline{h_w}|  \]
should be an $L^{q}_{loc}(\Omega')$ function.  That is essentially the same condition,  though notice that $h$ will satisfy the equation
\begin{eqnarray*}
  h_\wbar & = &  J(w,h) {\cal A}\Big(\frac{h_w}{J(w,h)}\Big).
  \end{eqnarray*}
It is not obvious how one unravels this equation to express $h_\wbar$ as a function of $h_w$ other than locally using the implicit function theorem,  but note that $ |h_\wbar |= J(w,h)| {\cal A}\Big(\frac{h_w}{J(w,h)}\Big)|\leq k |h_w|$ still holds.

\section{Fully nonlinear Beltrami systems.}

We now discuss the immediate consequences for the fully nonlinear Beltrami system.  Thus we consider the equation
\begin{equation}\label{eqn12}
f_\zbar={\cal H}(z,f,f_z)
\end{equation}  with ellipticity conditions that there is $k<1$ so  that for all $z\in \Omega$, $\zeta\in \Omega'$ and all $\zeta,\eta\in\IC$ 
\[|{\cal H}(z,w,\zeta)-{\cal H}(z,w,\eta)|\leq k |\zeta-\eta| \]
${\cal H},(z,\eta,0)\equiv 0$, with additional measurability condition ${\cal H}(z,w,\eta):\Omega\times\Omega'\times \IC\to\IC$ is Lusin measurable,  see \cite[Chapters 7 \& 8]{AIM} for the most general requirements.  Roughly,  the condition of Lusin measurability assures us that the function is measurable in each variable independently.  

The next result is clear from what is done above.

\begin{corollary}
Suppose that ${\cal H}(z,f,f_z):\Omega\times\Omega'\times \IC\to\IC$ has ellipticity bound $k$ and is linear at $\infty$ on the solution $f$.  That is there are constants $a,b\in\IC$ such that $|a|+|b|<1$ and 
\begin{equation}
|{\cal H}(z,f,f_z) - a f_z -b \overline{ f_z} | \in L^{q}_{loc}(\Omega), \hskip10pt \mbox{for some $q>1+k$}.
\end{equation}
Then every $W^{1,1}_{loc}(\Omega)$ solution $f$ to (\ref{eqn12}) has in fact $f\in W^{1,s}_{loc}(\Omega)$ for all $s<1+1/k$.  Hence $f$ is quasiregular.
\end{corollary}
\noindent{\bf Proof.} We may suppose $q<2$.  Put $u(z)={\cal H}(z,f,f_z) - a f_z -b \overline{ f_z}$.  Our hypotheses put $u\in L^{q}_{loc}(\Omega)$ and we observe $f_\zbar = af_z+b\overline{f_z}+u(z)$.  The proof of Theorem \ref{theorem3} (without using induction) puts $f\in W^{1,q}_{loc}(\Omega)$.  Then Astala's theorem \cite{Ast} puts $f\in W^{1,s}_{loc}(\Omega)$ and the ellipticity estimate gives us this for all $s<1+1/k$.  Observe that $1+1/k>2$. \hfill $\Box$

\medskip

However using an induction to improve regularity offers other alternative results for the nonlinear case.  We formulate one such now.  Suppose that the ellipticity and measurability conditions following (\ref{eqn12}) hold.  Specifically:   

\begin{enumerate}
\item There is $k<1$ so  that for all $z\in \Omega$, $\zeta\in \Omega'$ and all $\zeta,\eta\in\IC$ 
\[|{\cal H}(z,w,\zeta)-{\cal H}(z,w,\eta)|\leq k |\zeta-\eta| \]
\item ${\cal H},(z,\eta,0)\equiv 0$,
\item ${\cal H}(z,w,\eta):\Omega\times\Omega'\times \IC\to\IC$ is Lusin measurable,
\item That there are $a,b\in \IC$,  $|a|+|b|\leq k$ and
\begin{equation}
{\cal H}(z,\eta,\zeta) = a\zeta+b\overline{\zeta} + U(z,\eta,\zeta)
\end{equation}
and there is $\alpha<1$ and $u\in L^{q}_{loc}(\Omega)$ for some $q>1+k$ such that
\begin{equation}\label{eqn15}
|U(z,\eta,\zeta)| \leq A |\zeta|^\alpha + B |\eta|^{2\alpha} + u(z). 
\end{equation}
\end{enumerate}

\begin{theorem} Let $f:\Omega\to\Omega'$ be a $W^{1,1}_{loc}(\Omega)$ solution to the nonlinear equation(\ref{eqn12}) satisfying (1)-(4) above.  Then $f\in W^{1,2}_{loc}(\Omega)$ is a $\frac{1+k}{1-k}$-quasiregular mapping.
\end{theorem}

The proof is basically the same.  Initially we have $f\in L^{q}_{loc}(\Omega)$ for all $q<2$ by Sobolev embedding.  So the first two terms will self improve in an induction. The last term sets the limit of the improvability of the Sobolev exponent,  but with the ellipticity estimate and Astala's theorem this will be enough to guarantee quasiregularity.  There are other ways to reformulate the bounds on condition(\ref{eqn15}) so as to achieve a similar result.  We leave these for the reader to consider.

\bigskip 
 \noindent G. Martin, Massey University,   New Zealand,  g.j.martin@massey.ac.nz
\end{document}